\input amstex 
\input amsppt.sty
\magnification=\magstep1
\def\nmb#1#2{#2}         
\def\idx{}               
\def\ign#1{}             

\redefine\o{\circ}

\define\de{\delta}

\define\ph{\varphi}

\define\ps{\psi}
\define\om{\omega}

\define\De{\Delta}

\define\Om{\Omega}
\predefine\ii\i
\redefine\i{^{-1}}
\define\row#1#2#3{#1_{#2},\ldots,#1_{#3}}

\define\Der{\operatorname{Der}}
\define\Hom{\operatorname{Hom}}

\redefine\L{{\Cal L}}
\def\today{\ifcase\month\or
 January\or February\or March\or April\or May\or June\or
 July\or August\or September\or October\or November\or December\fi
 \space\number\day, \number\year}
\hyphenation{ho-mo-mor-phism}
\topmatter
\title  
Towards the Chern-Weil Homomorphism  \\
in Non Commutative Differential Geometry \endtitle
\author  Andreas Cap\\
Peter W. Michor 
\endauthor
\leftheadtext{\smc Cap, Michor}
\rightheadtext{\smc Non commutative Chern-Weil homomorphism}
\thanks{Supported by Project P 7724 PHY 
of `Fonds zur F\"orderung der wissenschaftlichen 
Forschung'\hfill}\endthanks
\affil
Institut f\"ur Mathematik, Universit\"at Wien,\\
Strudlhofgasse 4, A-1090 Wien, Austria.
\endaffil
\address{Institut f\"ur Mathematik, Universit\"at Wien,
Strudlhofgasse 4, A-1090 Wien, Austria.}\endaddress
\email {cap\@awirap.bitnet, cap\@pap.univie.ac.at, 
michor\@awirap.bitnet, \newline michor\@pap.univie.ac.at} \endemail


\def\rightheadtext{\smc Non commutative Chern-Weil homomorphism}
\document

\head \nmb0{1}. Introduction \endhead
In this short review article we sketch some developments which should 
ultimately lead to the analogy of the Chern-Weil homomorphism for 
principle bundles in the realm of non commutative differential 
geometry. Principal bundles there should have Hopf algebras as 
structure `cogroups'. Since the usual machinery of Lie algebras, 
connection forms, etc\., just is not available in this setting, we 
base our approach on 
the Fr\"olicher--Nijenhuis bracket. See \cite{9} for an account of 
the classical theory using this approach.

In this paper we give an outline of the construction of a non
commutative analogy of the Fr\"olicher--Nijenhuis bracket as well as
some simple applications. For simplicity we work in a purely
algebraic setting but the whole theory can also be developed for
topological algebras as well as for the so called convenient algebras
(see \cite{5}) which are best suited for
differentiation and take care of completed tensor products. 
For a detailed exposition in the latter setting see \cite{1} and 
\cite{2}.

\head\nmb0{2}. Universal differential forms \endhead
\subhead \nmb.{2.1}\endsubhead
Let $A$ be a unital associative algebra over a commutative field $K$ of 
characteristic zero. Then the graded algebra $\Om _*(A)$ of universal 
differential forms over $A$ is constructed as follows (see
\cite{6} and \cite{7}): The tensor product $A\otimes A$ is an 
$A$--bimodule and the multiplication map $\mu :A\otimes A\to A$ is a
bimodule homomorphism by associativity. Hence the kernel of $\mu$ is
an $A$--bimodule which we denote by $\Om _1(A)$. We define
$d:A\to \Om _1(A)$ by $d(a)=1\otimes a-a\otimes 1$. 

The map $d$ can be characterized by a universal property as follows:
Let $M$ be an $A$--bimodule. A linear map $D:A\to M$ is called a
derivation if and only if for any $a,b\in A$ we have 
$D(ab)=D(a)\cdot b+a\cdot D(b)$. 
Obviously $d$ is a derivation and thus for any bimodule
homomorphism $\ph :\Om _1(A)\to M$ the map $\ph \o d:A\to M$ is a
derivation. It can be proved that any derivation is of this form:

\proclaim{\nmb.{2.2}. Proposition} For any $A$-bimodule $M$ the canonical 
linear mapping 
$d^*:\Hom^A_A(\Om_1(A),M)\to \Der(A;M)$, given by 
$\ph\mapsto \ph\o d$ is an isomorphism.
\endproclaim

Clearly the module $\Om _1(A)$  is determined by this universal
property up to canonical isomorphism.

\subhead \nmb.{2.3}\endsubhead
Now we define the spaces of differential forms of higher degree by 
$\Om _k(A):=\Om _1(A)\otimes _A\dots\otimes _A\Om _1(A)$ ($k$ factors).
Moreover we put $\Om _0(A)=A$ and $\Om (A)=\oplus_{k=0}^\infty\Om
_k(A)$. Then $\Om (A)$ is a graded algebra with the tensor product as
multiplication. Next put $\bar A:=A/K$, where the ground field $K$ is
identified with the multiples of the unit of $A$. Then one proves
that the map $a\otimes \bar b\mapsto ad(b)$ is an isomorphism between
$A\otimes \bar A$ and $\Om _1(A)$. Consequently the map $a_0\otimes
\bar a_1\otimes\dots\otimes \bar a_k\mapsto a_0d(a_1)\dots d(a_k)$
defines an isomorphism $A\otimes \oversetbrace \text{$k$-times} 
\to{\bar A\otimes \cdots \otimes \bar A} \to \Om_k(A)$, and it turns
out that $a_0d(a_1)\dots d(a_k)\mapsto d(a_0)d(a_1)\dots d(a_k)$
gives a well defined map $d:\Om _k(A)\to \Om _{k+1}(A)$ for any $k$.
Then it can be shown that $(\Om (A),d)$ is a graded differential
algebra, i\.e\. that $d(\om_p \om_q) = d\om_p \om_q + (-1)^p\om_p d\om_q$
for all $\om_p\in\Om_p(A)$ and $\om_q\in\Om_q(A)$. Again this algebra
is characterized by a universal property:

\proclaim{\nmb.{2.4}. Proposition}
Let $(B=\oplus _{k=0}^\infty B_k,\de )$ be an arbitrary unital graded
differential algebra, $\ph _0:A\to B_0$ a homomorphism of unital
algebras. Then there is a unique homomorphism $\ph :\Om (A)\to B$ of
graded differential algebras which restricts to $\ph _0$ in degree
zero. 
\endproclaim

In particular this result shows that the construction of the algebra
of universal differential forms is functorial.

\head\nmb0{3}. Construction of the Fr\"olicher--Nijenhuis
bracket\endhead
The construction is based on the classification of all graded
derivations of the graded algebra $\Om (A)$.
\subhead\nmb.{3.1} Definition\endsubhead
The space $\operatorname{Der}_k\Om(A)$ consists of
all \idx{\it (graded) derivations} of degree $k$, i.e.
all bounded linear mappings 
$D:\Om(A) \to  \Om(A)$ with $D(\Om_\ell(A)) \subset
\Om_{k+\ell}(A)$ and 
$D(\ph \ps) = D(\ph) \ps +(-1)^{k\ell}\ph
D(\ps)$ for $\ph \in \Om_\ell(A)$.

\proclaim{Lemma} The space 
$\Der\Om(A) = \bigoplus_k\Der_k\Om(A)$ is a
graded Lie algebra with the graded commutator 
$[D_1,D_2] := D_1\o D_2 - (-1)^{k_1k_2}D_2\o D_1$ as bracket.
So the bracket is graded anticommutative,
$[D_1,D_2] = -(-1)^{k_1k_2}[D_2,D_1]$, and satisfies the graded
Jacobi identity
$[D_1,[D_2,D_3]] = [[D_1,D_2],D_3] + (-1)^{k_1k_2}[D_2,[D_1,D_3]]$.
\endproclaim

\subheading{\nmb.{3.2}} A derivation $D \in \Der_k\Om(A)$ is called
\idx{\it algebraic} if $D\mid \Om_0(A) = 0$.
Then $D(a\om) =
aD(\om)$ and $D(\om a)=D(\om)a$ for $a\in A$, so $D$ restricts to a 
bounded bimodule homomorphism, an element of 
$\Hom^A_A(\Om_l(A),\Om_{l+k}(A))$.
Since we have
$\Om_l(A) = \Om_1(A) \otimes_A\dots\otimes_A\Om_1(A)$
and since for a product of one forms we have
$D(\om_1\otimes _A\dots\otimes _A\om_l)=\sum_{i=1}^l(-1)^{ik}
\om_1\otimes _A\dots\otimes _A D(\om_i)\otimes _A\dots\otimes_A\om_l$,
the derivation $D$ is uniquely determined by its restriction
$K:=D|\Om_1(A)\in\Hom^A_A(\Om_1(A),\Om_{k+1}(A))$;
we write $D=j(K)=j_K$ to express this dependence. Note the defining 
equation $j_K(\om)= K(\om)$ for $\om\in\Om_1(A)$.
Since it will be very important in the sequel we will use the 
notation
$\Om^1_k=\Om^1_k(A):=\Hom^A_A(\Om_1(A),\Om_k(A))$ and 
$\Om^1_*=\Om^1_*(A)=\bigoplus_{k=0}^\infty\Om^1_k(A)$.

It can be shown that for any $K\in \Om ^1_k(A)$ the formula 
$$j_K(\om_0\otimes _A  \dots\otimes _A  \om_\ell)
     =\sum_{i=0}^\ell(-1)^{ik}\om_0\otimes _A  \dots
     \otimes _A  K(\om_i)\otimes _A  \dots \otimes _A \om_k$$
for $\om _i\in \Om_1(A)$
defines an algebraic graded derivation $j_K \in \Der_k\Om(A)$ and any
algebraic derivation is of this form. Thus $K\mapsto j_K$ is an
isomorphism from $\Om ^1_*(A)$ to the space of algebraic graded
derivations of $\Om (A)$ Since the graded commutator of two algebraic
derivations is clearly again algebraic we can define a graded Lie
bracket $[\quad,\quad]^{\De}$ on the space $\Om ^1_*(A)$ by  
$j([K,L]^{\De}):= [j_K,j_L]$. This bracket is called the \idx{\it 
algebraic bracket}; it is an analogy of the one used in \cite{3}.

\subhead\nmb.{3.3} \endsubhead
The differential $d$ is a graded derivation of $\Om (A)$ of degree
one which is not algebraic.
In analogy to the well known formula for the Lie derivative along
vector fields we now define the {\it Lie derivative\/} along a field
valued form $K\in \Om ^1_k(A)$ by $\L_K := [j_K,d]\in \Der_k\Om(A)$.
Then one proves that for any derivation $D \in \Der_k\Om(A)$ there
are unique elements $K \in \Om^1_k$ and $L \in \Om^1_{k+1}$ such that 
$D = \L_K + j_L$. Moreover $L=0$ if and only if $[D,d]=0$ and $D$ is 
algebraic if and only if $K=0$.

For elements $K \in \Om^1_k$ and $L \in \Om^1_\ell$ one immediately
verifies that $[[\L_K,\L_L],d] =0$, so we have
$[\L(K),\L(L)] = \L([K,L])$
for a uniquely defined 
$[K,L] \in \Om^1_{k+\ell}$. This
vector valued form $[K,L]$ is called the \idx{\it abstract
Fr\"olicher-Nijenhuis bracket} of $K$ and $L$. Clearly this bracket
defines a graded Lie algebra structure on the space $\Om ^1_*(A)$. 

\head\nmb0{4}. Distributions and Integrability\endhead
\subheading{\nmb.{4.1}. Distributions} By a \idx{\it distribution} 
in an algebra $A$ we mean a sub-$A$-bimodule $\Cal D$ of $\Om_1(A)$.

The distribution $\Cal D$ is called \idx{\it globally integrable} if 
there exists a sub algebra $B$ of $A$ such that
$\Cal D$ is the subspace generated by $A(d(B))$ and $d(B)A$.

The distribution $\Cal D$ is called \idx{\it splitting} if it is a 
direct summand in $\Om _1(A)$ or equivalently if there is a 
projection $P\in\Om^1_1(A)=\Hom^A_A(\Om_1(A),\Om_1(A))$
onto $\Cal D$, i.e. $P\o P=P$ and $\Cal D = P(\Om_1(A))$.  Then there 
is a complementary sub module $\ker P\subset \Om_1(A)$.

The distribution $\Cal D$ is called \idx{\it involutive} if the 
ideal $(\Cal D)_{\Om_*(A)}$ generated by $\Cal D$ in the
graded algebra $\Om_*(A)$ is stable under $d$, i\.e\. if
$d(\Cal D)\subset (\Cal D)_{\Om_*(A)}$.

\subheading{\nmb.{4.2}. Comments}
One should think of this as follows: In ordinary differential geometry 
$\Cal D$ should be the $A$-bimodule of those 1-forms which annihilate 
the sub bundle of $TM$. Global integrability then means that it is 
integrable and that the space of functions which are constant along 
the leaves of the foliation generates those forms. This is a strong 
condition: There are foliations where this space of functions 
consists only of the constants, and this can be embedded into any 
manifold. So in $C^\infty(M)$ there are always involutive 
distributions which are not globally integrable. To prove some
Frobenius theorem a notion of local integrability would be necessary.

\subheading{\nmb.{4.3} Curvature and cocurvature}  
Let $P\in\Om^1_1(A)=\Hom^A_A(\Om_1(A),\Om_1(A))$
be a projection, then the image $P(\Om_1(A))$ is a splitting distribution, 
called the \idx{\it vertical distribution} of $P$ and 
the complement $\ker P$ is also a splitting distribution, called the 
\idx{\it horizontal\ign{ distribution}} one.
$\bar P:= Id_{\Om_1(A)}- P$ is a projection onto  the 
horizontal distribution.

We consider now the Fr\"olicher-Nijenhuis bracket $[P,P]$ of $P$ and 
define 
$$\alignat2
R &= R_P = [P,P]\o P &\qquad &\text{ the \it curvature.}\\
\bar R &= \bar R_P = [P,P]\o \bar P &\qquad &\text{ the \it cocurvature,}
\endalignat$$
The \idx{curvature} and the \idx{cocurvature} are elements of   
$\Om_2^1(A)=\Hom^A_A(\Om_1(A),\Om_2(A))$. The curvature kills 
elements of the horizontal distribution, so it is \idx{\it vertical}.
The cocurvature kills elements of the vertical distribution.

Then one proves:
\proclaim{ Proposition} With $P$, $R$ and $\bar R$ as above we have:

1. (Bianchi Identity)
$$\align &[P ,R+\bar R] = 0 \\
         &2[R,P ] = j_R\bar R + j_{\bar R}R,
\endalign $$
where the insertion operators are extended to $\Om ^1_*(A)$ in the 
obvious way.

2. The curvature $R$ is zero if and only if the horizontal 
distribution is involutive. The cocurvature $\bar R$ 
is zero if and only if the vertical distribution
$P(\Om_1(A))$ is involutive.  
\endproclaim

\enddocument